\documentclass[12pt]{article}

\usepackage{amsmath,amssymb,enumerate,bbm}
\newcommand{\tmop}[1]{\ensuremath{\operatorname{#1}}}
\newenvironment{itemizedot}{\begin{itemize} }{\end{itemize}}

\newtheorem{theorem}{Theorem}[section]

\newtheorem{lemma}[theorem]{Lemma}

\numberwithin{equation}{section}
\newenvironment{proof}{\noindent\textbf{Proof\ }}{\hspace*{\fill}$\Box$\medskip}

 \makeatletter
 
 \newcommand{\Rmnum}[1]{\expandafter\@slowromancap\romannumeral #1@}
 \makeatother
\begin{document}

\title{On kaleidoscopic pseudo-randomness of finite Euclidean graphs}\author{Le Anh
Vinh\\
Mathematics Department\\
Harvard University\\
Cambridge, MA 02138, US}\maketitle

\begin{center}
\end{center}

\begin{abstract}
  D. Hart, A. Iosevich, D. Koh, S. Senger and I. Uriarte-Tuero (2008) showed that the distance graphs has kaleidoscopic pseudo-random property, i.e. sufficiently large subsets of $d$-dimensional vector spaces over finite   fields contain every possible finite configurations.  In this paper we study the kaleidoscopic pseudo-randomness of finite Euclidean graphs using probabilistic methods.  
\end{abstract}

\section{Introduction}

Let $\mathbbm{F}_q$ denote the finite field with $q$ elements where $q \gg 1$
is an odd prime power. For a fixed $a \in \mathbbm{F}_q^{\ast}$, the distance graph $G_q(a)$ (also known as finite Euclidean graphs in \cite{medrano}) in $\mathbbm{F}_q^d$ is defined as the graph with vertex
set $\mathbbm{F}_q^d$ and the edge set
\[ E =\{(x, y) \in \mathbbm{F}_q^d \times \mathbbm{F}_q^d \mid x \neq y, ||x - y|| = a\}, \]
where $||.||$ is the analogue of Euclidean distance $||x|| = x_1^2+\ldots+x_d^2$.

In other words, consider the set of colors $L =\{c_1, \ldots, c_{q - 1} \}$
corresponding to elements of $\mathbbm{F}_q^{\ast}$. We color the complete
graph $K_{q^d}$ with vertex set $\mathbbm{F}_q^d$ by $q - 1$ colors such that
$(x, y) \in \mathbbm{F}_q^d \times \mathbbm{F}_q^d$ is colored by the color
$c_i$ if $||x - y|| = i$. Denote this resulting family of graphs, with respect to the
above coloring, by $G_q^{\Delta}$ where $q$ runs over powers of odd primes.

In \cite{medrano}, Medrano et al. studied the spectrum of these graphs and showed that these graphs are asymptotically Ramanujan graphs. In \cite{bannai}, Bannai, Shimabukuro and Tanaka showed that the distance graphs over finite fields are always asymptotically Ramanujan for a more general setting (i.e. they replace the Euclidean distance above by nondegenerated quadratic forms). The author recently applied these results to several interesting combinatorial problems, for example to tough Ramsey graphs (with P. Dung) \cite{vinh-dung}, to the Erd\"os distance problem \cite{vinh}, to integral graphs (with Si Li) \cite{si-vinh}, also to a Szemeredi-Trotter type theorem and to a sum-product estimate \cite{vinh-sp} (see \cite{vu} for related results). 

The definition of kaleidoscopic pseudo-randomness below follows Hart et al. \cite{hart et al}. We say that the family of graphs $\{G_j \}_{j = 1}^{\infty}$ with the sets of
colors $L_j =\{c_j^1, c_j^2, \ldots, c_j^{|L_j |} \}$ and the edge set
$\mathcal{E}_j = \bigcup_{i = 1}^{|L_j |} \mathcal{E}_j^i$, with
$\mathcal{E}_j^i$ corresponding to the color $c_j^i$, is kaleidoscopically
pseudo-random if there exist constants $C, C' > 0$ such that the following
conditions are satisfied:
\begin{itemizedot}
  \item
  \begin{equation}
    |G_j | \rightarrow \infty \tmop{as} j \rightarrow \infty .
  \end{equation}
  \item For any $j$ and $1 \leqslant i, i' \leqslant |L_j |$ then
  \begin{equation}
    \frac{1}{C'} |\mathcal{E}_j^{i'} | \leqslant |\mathcal{E}_j^i | \leqslant
    C' |\mathcal{E}_j^{i'} |.
  \end{equation}
  \item $G_j$ is asymptotically complete in the sense that
  \begin{equation}
    \lim_{j \rightarrow \infty} \frac{\binom{|G_j |}{2} - \sum_{i = 1}^{|L_j |}
    |\mathcal{E}_j^i |}{\binom{|G_j |}{2}} = 0.
  \end{equation}
  \item If $1 \leqslant k - 1 \leqslant n$ and $L_j' \subset L_j$, with $|L_j'
  | \leqslant |L_j | - \binom{k}{2} + n$, then any subgraph $H$ of $G_j$ of order
  \begin{equation}\label{14}
    \geqslant C|G_j |^{\frac{k - 1}{k}} |L|^{\frac{n}{k}},
  \end{equation}
  contains every possible subgraph with $k$ vertices and $n$ edges with an
  arbitrary edge color distribution from $L_j'$.
\end{itemizedot}

In \cite{hart et al}, Hart et al. systematically studied various properties of
the distance graph. Using Fourier analysis, they proved the following result (see \cite{hart et
al} for the motivation and applications of this result).

\begin{theorem} \label{main 1} (\cite{hart et al})
  The family of graph $\{G_q^{\Delta} (d)\}$ is kaleidoscopically
  pseudo-random.
\end{theorem}

The first three items in the definition of kaleidoscopically pseudo-randomness above are easy special cases of the following lemma (which is implicit in several papers, for example, see \cite{hart et al, iosevich, medrano, vinh}). 

\begin{lemma} (Lemma 1.2 in \cite{hart et al})
For any $t \in \mathbbm{F}_q$, then
\begin{equation}
|\{(x,y) \in \mathbbm{F}_q^d \times \mathbbm{F}_q^d : ||x-y|| = t \}| = \left\{
        \begin{array}{ll}
        (2+o(1))q^{2d-1} & \mbox{if} \ d=2, t= 0\\
        (1+o(1))q^{2d-1} & \mbox{otherwise.}
    \end{array}  \right.
\end{equation}
where $o(1)$ means that the quantity goes to $0$ as $q$ goes to infinity. 
\end{lemma}

Hart et al. (\cite{hart et al}) derived the item (\ref{14}) from the following estimate (This estimate is just a graph theoretic translation of Theorem 1.4 in \cite{hart et al}).

\begin{theorem}\label{main 2}
Consider the color set $L =\{c_1, \ldots, c_{q - 1} \}$
coressponding to elements of $\mathbbm{F}_q^{\ast}$. Color the complete graph $K_{q^d}$ by the Euclidean distance as the above. Let $E \subset \mathbbm{F}_q^d$, $d \geq 2$. Suppose that $1 \leq k-1 \leq n \leq d$ and 
\begin{equation}
|E| \geq Cq^{d\frac{k-1}{k}}q^{\frac{n}{k}}
\end{equation}
with a sufficiently large constant $C > 0$. Then for any subgraph $H$ with $k$ vertices and $n$ edges with an arbitrary edge color distribution from $L' \subset L$ ($|L'| \leq |L| - \binom{k}{2}+n$), we have 
\begin{equation}
(1+o(1))|E|^kq^{-n}
\end{equation}
copies of $H$ (with vertex ordering) in $E$. 
\end{theorem}

In this paper, we will study more general families of distance graphs. Let $Q$ be a non-degenerate quadratic form on $\mathbbm{F}_q^d$. The finite Euclidean graph $E_q(n,Q,a)$ is defined as the graph with vertex set $\mathbbm{F}_q^d$ and the edge set
\begin{equation}
E = \{(x,y) \in \mathbbm{F}_q \times \mathbbm{F}_q \mid x \neq y, Q(x-y) = a \}.
\end{equation}

Similarly, consider the set of colors $L =\{c_1, \ldots, c_{q - 1} \}$
corresponding to elements of $\mathbbm{F}_q^{\ast}$. We color the complete
graph $K_{q^d}$ with vertex set $\mathbbm{F}_q^d$ by $q - 1$ colors such that
$(x, y) \in \mathbbm{F}_q^d \times \mathbbm{F}_q^d$ is colored by the color
$c_i$ if $Q(x - y) = i$. Denote this resulting family of graphs, with respect to this  coloring, by $G_q^{Q}$ where $q$ runs over powers of odd primes.

The main result of this paper is the following similar result of Theorem \ref{main 2} for finite Euclidean graphs.

\begin{theorem} \label{main}
Consider the color set $L =\{c_1, \ldots, c_{q - 1} \}$
coressponding to elements of $\mathbbm{F}_q^{\ast}$. Color the complete graph $K_{q^d}$ by the non-degenerate quadratic form $Q$ as the above. Let $E \subset \mathbbm{F}_q^d$, $d \geq 2$. Suppose that $1 \leq k-1 \leq n \leq d$ and 
\begin{equation}
|E| \gg q^{\frac{d-1}{2}+k-1}.
\end{equation}
Then for any subgraph $H$ with $k$ vertices and $n$ edges with an arbitrary edge color distribution from $L$, we have 
\begin{equation}
(1+o(1))|E|^kq^{-n}
\end{equation}
copies of $H$ (with vertex ordering) in $E$. 
\end{theorem}

The result is only non-trivial in the range $d \geq 2(k-1)$. Note that Theorem \ref{main} is stronger than Theorem \ref{main 2} in this range. Moreover, different from Hart et al. \cite{hart et al}, our proof uses probabilistic methods. The rest of this paper is organized as follows. In
Section 2, we establish a theorem about the number of small colored subgraphs
in pseudo-random coloring of a complete graph. Using this theorem, we give a proof of Theorem \ref{main} and also discuss
similar results in more general settings in the last section.

\section{Pseudo-random coloring}

We call a graph $G = (V, E)$ $(n, d, \lambda)$-graph if $G$ is a $d$-regular
graph on vertices with the absolute values of each of its eigenvalues but the
largest one is at most $\lambda$. It is well-known that if $\lambda \ll d$
then a $(n, d, \lambda)$ graph behaves similarly as a random graph $G_{n, d /
n}$. Precisely, we have the following result (see Corollary 9.2.5 in
\cite{alon-spencer}).

\begin{theorem} \label{tool 1}
  (\cite{alon-spencer}) Let $G$ be an $(n, d, \lambda)$-graph. For every set
  of vertices $B$ and $C$ of $G$, we have
  \begin{equation}
    |e (B, C) - \frac{d}{n} |B||C\| \leqslant \lambda \sqrt{|B||C|},
  \end{equation}
  where $e (B, C)$ is the number of edges in the induced bipartite subgraph of
  $G$ on $(B, C)$ (i.e. the number of ordered pairs $(u, v)$ where $u \in B$,
  $v \in C$ and $u v$ is an edge of $G$).
\end{theorem}

Let $H$ be a fixed graph of order $s$ with $r$ edges and with automorphism
group $\tmop{Aut} (H)$. It is well-known that for every constant $p$ the
random graph $G (n, p)$ contains
\begin{equation}
  (1 + o (1)) p^r (1 - p)^{(^s_2) - r} \frac{n^s}{| \tmop{Aut} (H) |}
\end{equation}
induced copies of $H$. Alon extended this result to $(n, d, \lambda)$-graph.
He proved that every large subset of the set of vertices of a $(n, d,
\lambda)$-graph contains the ``correct'' number of copies of any fixed small
subgraph (Theorem 4.10 in \cite{krivelevich-sudakov}).

\begin{theorem}\label{tool 2}
  (\cite{krivelevich-sudakov}) Let $H$ be a fixed graph with $r$ edges, $s$
  vertices and maximum degree $\Delta$, and let $G = (V, E)$ be an $(n, d,
  \lambda)$-graph, where, say, $d \leqslant 0.9 n$. Let $m < n$ satisfies $m
  \gg \lambda \left( \frac{n}{d} \right)^{\Delta}$. Then, for every subset $U
  \subset V$ of cardinality $m$, the number of (not necessrily induced) copies
  of $H$ in $U$ is
  \begin{equation}
    (1 + o (1)) \frac{m^s}{| \tmop{Aut} (H) |} \left( \frac{d}{n} \right)^r .
  \end{equation}
\end{theorem}

Suppose that a graph $G$ of order $n$ is colored by $t$ colors. Let $G_i$ be
the induced subgraph of $G$ on the $i^{\tmop{th}}$ color. We call a
$t$-colored graph $G$ $(n, d, \lambda)$-r.c (regularly colored)
graph if $G_i$ is a $(n, d, \lambda)$-regular graph for each color $i \in \{1,
\ldots, t\}$. We present here a variant of Theorem \ref{tool 2} that very large subset of
the vertex set of a $(n, d, \lambda)$-r.c graph contains the ``correct'' number of
copies of any fixed small colored graph.

\begin{theorem} \label{tool 3}
  Let $H$ be a fixed $t$-colored graph with $r$ edges, $s$ vertices, maximum
  degree $\Delta$ with automorphism group (with respect to coloring)
  $\tmop{Aut}_c (H)$, and let $G$ be a $t$-colored graph of order $n$.
  Suppose that $G$ is an $(n, d, \lambda$)-r.c graph, where, say, $d \ll n$.
  Let $m < n$ satisfies $m \gg \lambda \left( \frac{n}{d} \right)^{\Delta}$.
  Then, for every subset $U \subset V$ of cardinality $m$, the number of (not
  necessrily induced) copies of $H$ in $U$ is
  \begin{equation}
    (1 + o (1)) \frac{m^s}{| \tmop{Aut}_c (H) |} \left( \frac{d}{n} \right)^r.
  \end{equation}
  If we take the ordering of vertex set into account then the number of copies of $H$ in $U$ is
  \begin{equation}
    (1 + o (1)) m^s \left( \frac{d}{n} \right)^r.
  \end{equation}
\end{theorem}

The proof of this theorem is similar to the proof of Theorem 4.10 in
\cite{krivelevich-sudakov}. We give a detail proof here for the sake of
completeness.

\begin{proof}
  To prove the theorem, consider a random one-to-one mapping of the set of
  vertices of $H$ into the set of vertices $U$. Let $M (H)$ denote the event
  that every edge of $H$ is mapped on an edge of $K_n$ with corresponding
  color. We say that the mapping is an embedding of $H$ in such a case. It
  suffices to prove that
  \begin{equation} \label{5}
    \Pr (M (H)) = (1 + o (1)) \left( \frac{d}{n} \right)^r .
  \end{equation}
  We prove (\ref{5}) by induction on the number of edges $r$. The base case $(r =
  0$) is trivial. Suppose that $(\ref{5})$ holds for all colored graphs with less
  than $r$ edges. Let $u v$ be an edge of $H$. Let $H_u, H_v, H_{\{u, v\}}$ be the induced
  subgraph of $H$ on the vertice set $V (H) -\{u\}$, $V (H) -\{v\}$, $V (H) -\{u, v\}$, and let $H_{u v}$ be the graph obtained
  from $H$ by removing the edge $u v$. We have
  \begin{equation}
    \text{$\Pr (M (H_{u v})) = \Pr$} (M (H_{u v}) |M (H_{\{u, v\}})) . \Pr (M
    (H_{\{u, v\}})) .
  \end{equation}
  Let $r_1$ be the number of edges of $H_{\{u, v\}}$. Since (\ref{5}) holds for
  $H_{u v}$ and $H_{\{u, v\}}$, we have $\Pr (M (H_{u v})) = (1 + o (1))
  \left( \frac{d}{n} \right)^{r - 1}$ and $\Pr (M (H_{\{u, v\}}) = (1 + o (1))
  \left( \frac{d}{n} \right)^{r_1}$. Thus, we have
  \begin{equation}
    \text{$\Pr$} (M (H_{u v}) |M (H_{\{u, v\}})) = (1 + o (1)) \left(
    \frac{d}{n} \right)^{r - r_1 - 1} .
  \end{equation}
  For an embedding $f_1$ of $H_{\{u, v\}}$ in $U$, let $\phi (u, f_1)$,
  $\phi (v, f_1)$ and $\phi (u v, f_1)$ be the number of extensions of $f_1$ to an embedding of
  $H_u$, $H_v$ and $H_{u v}$ in $U$, respectively. Note that an
  extension $f_u$ of $f_1$ to an embedding of $H_u$ and an extension $f_v$ of
  $f_1$ to an embedding of $H_v$ give us a unique extension of $f_1$ to an
  embedding of $H_{u v}$ except that the image of $u$ in $f_u$ is the same as
  the image of $v$ in $f_v$. Thus, we have
  \begin{equation}
    \phi (u, f_1) \phi (v, f_1) - \min (\phi (u, f_1), \phi (v, f_1))
    \leqslant \phi (u v, f_1) \leqslant \phi (u, f_1) \phi (v, f_1) .
  \end{equation}
  Averaging over all possible extensions of $f_1$ to a mapping from $H_{u v}$
  to $U$, we have
  \begin{equation*}
    \frac{\phi (u, f_1) \phi (v, f_1) - \min (\phi (u, f_1), \phi (v,
    f_1))}{(m - s + 2) (m - s + 1)} \leqslant \Pr (M (H_{u v}) |f_1) \leqslant
    \frac{\phi (u, f_1) \phi (v, f_1)}{(m - s + 2) (m - s + 1)} .
  \end{equation*}
  Taking expectation over all embedding $f_1$, the middle term becomes
  \[\text{$\Pr$} (M (H_{u v}) |M (H_{\{u, v\}})) = (1 + o (1)) \left(
  \frac{d}{n} \right)^{r - r_1 - 1}.\]
  Note that $\min (\phi (u, f_1), \phi (v,
  f_1)) \leqslant m$ so we get
  \[ E_{f_1} (\phi (u, f_1) \phi (v, f_1) |M (H_{\{u, v\}})) = (1 + o (1))
     m^2 \left( \frac{d}{n} \right)^{r - r_1 - 1} + \delta, \]
  where $| \delta | \leqslant m$. We have $r - r_1 \leqslant 2 (\Delta - 1) +
  1$ and $m \gg \lambda \left( \frac{n}{d} \right)^{\Delta} = \Omega (
  \sqrt{d}) \left( \frac{n}{d} \right)^{\Delta}$, thus
  \begin{equation}
    m^2 \left( \frac{d}{n} \right)^{r - r_1 - 1} \geqslant m^2 \left(
    \frac{d}{n} \right)^{2 (\Delta - 1)} \gg \Omega (d) \left( \frac{n}{d}
    \right)^2 = \Omega (n^2 / d) \gg n > m.
  \end{equation}
  So $\delta$ is negligible and we get
  \begin{equation}
    E_{f_1} (\phi (u, f_1) \phi (v, f_1) |M (H_{\{u, v\}})) = (1 + o (1)) m^2
    \left( \frac{d}{n} \right)^{r - r_1 - 1} .
  \end{equation}
  Now, let $f$ be a random one-to-one mapping of $V (H)$ into $U$. Let $f_1$
  be a fixed embedding of $H_{\{u, v\}}$. Let $B \tmop{and} C$ be the set of
  all possible images of $u$ and $v$ over all possible extensions of $f_1$ to
  an embedding of $H_u$ and $H_v$ in $U$, respectively. Since $G$ is an $(n,
  d, \lambda)$-r.c graph, by Theorem \ref{tool 1}, the number of possible pairs $(u, v)$
  with $u \in B$ and $v \in C$ such that $u v$ is correctly colored is bounded
  by
  \begin{equation} \label{11}
    \frac{d}{n} \phi (u, f_1) \phi (v, f_1) \pm \lambda \sqrt{\phi (u, f_1)
    \phi (v, f_1)} .
  \end{equation}
  Thus, we have
  \[
    \Pr\ _f (M (H) |f_{|V (H) \backslash \{u, v\}} = f_1) = \frac{d}{n}
    \frac{\phi (u, f_1) \phi (v, f_1)}{(m - s + 2) (m - s + 1)} + \delta,
  \]
  where $| \delta | \leqslant \lambda \frac{\sqrt{\phi (u, f_1) \phi (v,
  f_1)}}{(m - s + 2) (m - s + 1)}$. Averaging over all possible embeddings
  $f_1$, we get
  \begin{eqnarray*}
    \Pr (M (H) |M (H')) & = & \frac{d}{n} \frac{E_{f_1} (\phi (u, f_1) \phi
    (v, f_1) |M (H_{\{u, v\}}))}{(m - s + 2) (m - s + 1)} + E_{f_1} (\delta)\\
    & = & (1 + o (1)) \left( \frac{d}{n} \right)^{r - r_1} + E_{f_1}
    (\delta),
  \end{eqnarray*}
  where the second lines follows from (\ref{11}). By Jensen's inequality, we have
  \begin{equation}
    |E_{f_1} (\delta) | \leqslant \lambda \frac{\sqrt{E (\phi (u, f_1) \phi
    (v, f_1))}}{(m - s + 2) (m - s + 1)} = (1 + o (1)) \frac{\lambda}{m}
    \left( \frac{d}{n} \right)^{(r - r_1 - 1) / 2},
  \end{equation}
  which is negligible to the first term (as $m \gg \lambda \left( \frac{n}{d}
  \right)^{\Delta} \geqslant \lambda \left( \frac{n}{d} \right)^{(r - r_1 + 1)
  / 2}$). Thus, we have
  \begin{equation}
    \Pr (M (H)) = \Pr (M (H) |M (H_{\{u, v\}})) \Pr (M (H_{\{u, v\}}) = (1 + o
    (1)) \left( \frac{d}{n} \right)^r .
  \end{equation}
  This completes the proof of the theorem.
\end{proof}

\section{Finite Euclidean graphs}
\subsection{Proof of Theorem \ref{main}}
In \cite{bannai}, Bannai, Shimabukuro and Tanaka studied the spectrum of distance graph $G_q(a)$ and showed that these graphs are asymptotically Ramanujan graphs. They proved the following result.

\begin{theorem} (\cite{bannai, medrano}) \label{tool 4}
Let $Q$ be a non-degenerate quadratic form on $\mathbbm{F}_q^d$. The finite Euclidean graph $E_q(d,Q,a)$ is regular of valency $(1+o(1))q^{d-1}$ for any $a \in \mathbbm{F}_q^{\ast}$. Let $\lambda$ be any eigenvalues of the graph $G_q(a)$ with $\lambda \neq$ valency of the graph then 
\begin{equation}
|\lambda| \leq 2q^{\frac{d-1}{2}}.
\end{equation}
\end{theorem}

From Theorem \ref{tool 4}, if we color the complete graph $G = K_{q^d}$ with vertex set $\mathbbm{F}_q^{d}$ by a non-degenerate quadratic form $Q$ as in Section 1 then the colored graph $G$ is a $(q^d, (1+o(1))q^{d-1},2q^{\frac{d-1}{2}})$-r.c. graph. The Theorem \ref{main} now follows immediately from Theorem \ref{tool 3}. 

\subsection{General distances}

The proof in \cite{hart et al} shows that the conclusion of Theorem \ref{main 2} holds with the Euclidean norm $||.||$ is replaced by any function $F$ with the property that the Fourier transform satisfies the decay estimates
\begin{equation}\label{cd1}
\left| \hat{F}_t(m) \right| = \left| q^{-d} \sum_{ x\in \mathbbm{F}_q^d : F(x) = t} \chi (-x . m) \right| \leqslant C q^{-(d+1)/2}
\end{equation}
and
\begin{equation} \label{cd2}
\left| \hat{F}_t(0,\ldots,0) \right| = \left| q^{-d} \sum_{ x\in \mathbbm{F}_q^d : F(x) = t} \chi (-x . (0,\ldots,0)) \right| \leqslant C q^{-1},
\end{equation}
where $\chi(s) = e^{2\pi i \text{Tr}(s)/q}$ and $m \neq (0,\ldots,0) \in \mathbbm{F}_q^d$ (recall that for $y \in \mathbbm{F}_q$, where $q=p^r$ with $p$ prime, the trace of $y$ is defined as $\text{Tr}(y) = y + y^p + \ldots + y^{p^{r-1}} \in \mathbbm{F}_q$).

Now we define the $F$-distance graph $G_F(q,d,j)$ with the vertex set $V = \mathbbm{F}_q^d$ and the edge set
\[E=\{(x,y) \in V \times V | x \neq y, F(x-y) = j\}.\]

Then the exponentials (or characters of the additive group $\mathbbm{F}_q^d$)
\begin{equation}
e_m(x) = \exp \left( \frac{2\pi i \text{Tr} (x.m)}{p} \right), 
\end{equation}
for $x, m \in \mathbbm{F}_q^d$, are eigenfunctions of the adjacency operator for the $F$-distance graph $G_F(q,d,j)$ corresponding to the eigenvalue
\begin{equation}\label{e}
\lambda_m = \sum_{F(x) = j} e_m(x) = q^d \hat{F}_j(-m).
\end{equation}
Thus, the decay estimates (\ref{cd1}) and (\ref{cd2}) are equivalent to
\begin{equation}\label{b1}
\lambda_m \leqslant C q^{(d-1)/2},
\end{equation}
for $m\neq (0,\ldots,0) \in \mathbbm{F}_q^d$, and
\begin{equation}\label{b2}
\lambda_{(0,\ldots,0)} \leqslant C q^{d-1}.
\end{equation}

Therefore, we can apply Theorem \ref{tool 4} to obtain similar results. 

\section*{Acknowledgements}

The research is performed during the author's visit at the Erwin Schr\"odinger International Institute for Mathematical Physics. The author would like to thank the ESI for hospitality and financial support during his visit.

\end{document}